\documentclass[11pt]{amsart}
\usepackage{tikz-cd}
\usepackage{amscd,amssymb,palatino,amsmath}

\usepackage[utf8]{inputenc}
\usepackage[english]{babel}

\usepackage{amssymb}
\usepackage{graphicx, color}
\usepackage[all,cmtip]{xy}
\usepackage{pb-diagram,pb-xy} 

\renewcommand{\H}{{\mathcal{H}}}

\newcommand{\Z}{\mathbb{Z}}

\newcommand{\R}{\mathbb{R}}

\newcommand{\g}{{\mathfrak g}}
\newcommand{\h}{{\mathfrak h}}

\newcommand{\frbd}[1]{\frac{\partial}{\partial #1}}

\newtheorem{proposition}{Proposition}
\newtheorem{theorem}[proposition]{Theorem}
\newtheorem{definition}[proposition]{Definition}
\newtheorem{lemma}[proposition]{Lemma}

\newtheorem{example}[proposition]{Example}

\begin{document}

\title{$b$-Structures on Lie groups and Poisson reduction}

\author{Roisin Braddell}\address{ Roisin Braddell,
Basque Center for Applied Mathematics,
Mazarredo, 14 E48009 Bilbao, Basque Country – Spain. \it{e-mail} rdempsey@bcamath.org
}
\author{Anna Kiesenhofer}\address{ Anna Kiesenhofer,
Department of Mathematics, EPFL, Lausanne \it{e-mail: anna.kiesenhofer@epfl.ch}
 }

\author{Eva Miranda}\address{ Eva Miranda,
Laboratory of Geometry and Dynamical Systems, Departament of Mathematics-IMTech, Universitat Polit\`{e}cnica de Catalunya, Avinguda del Doctor Mara\~{n}on 44-50, 08028 , Barcelona \& CRM Centre de Recerca Matem\`{a}tica, Campus de Bellaterra
Edifici C, 08193 Bellaterra, Barcelona
 \it{e-mail: eva.miranda@upc.edu
 }}\thanks{{ Eva Miranda is supported by the Catalan Institution for Research and Advanced Studies via an ICREA Academia Prize 2016, by the project PID2019-103849GB-I00 of MCIN/ AEI /10.13039/501100011033, the project AGAUR grant 2017SGR932 and by the Spanish State
Research Agency, through the Severo Ochoa and Mar\'{\i}a de Maeztu Program for Centers and Units
of Excellence in R\&D (project CEX2020-001084-M). Roisin Braddell is supported by the Severo Ochoa Program SEV-2017-0718, Basque Government BERC Program 2018-2021 and was supported by a predoctoral grant from UPC with the ICREA Academia project of Eva Miranda when this project was initiated. We acknowledge support from the Fondation Sciences Mathématiques de Paris via the Chaire d'Excellence of Eva Miranda supported by a public grant overseen by the French National Research Agency (ANR) as part of the "Investissements d'Avenir" program (reference: ANR-10-LABX-0098) to finance a research stay of Roisin Braddell and a research visit of Anna Kiesenhofer in Paris to start this project.}}

\begin{abstract} Motivated by the group of Galilean transformations and the subgroup of Galilean transformations which fix time zero, we introduce the notion of a $b$-Lie group as a pair $(G,H)$ where $G$ is a Lie group and $H$ is a codimension-one Lie subgroup. Such a notion allows us to give a theoretical framework for transformations of space-time where the initial time can be seen as a boundary. In this theoretical framework, we develop the basics of the theory and study the associated canonical $b$-symplectic structure on the $b$-cotangent bundle $^b {T}^\ast G$ together with its reduction theory. Namely, we extend the minimal coupling procedure to $^bT^*G/H$ and prove that the Poisson reduction under the cotangent lifted action of $H$ by left translations can be described in terms of the Lie Poisson structure on $\mathfrak{h}^\ast$ (where $\mathfrak{h}$ is the Lie algebra of $H$) and the canonical $b$-symplectic structure on
$^b {T}^\ast(G/H)$, where $G/H$ is viewed as a one-dimensional $b$-manifold having as critical hypersurface (in the sense of $b$-manifolds) the identity element.
\end{abstract}

\maketitle


\section{Introduction}

The Galilean group $G$ is the group of transformations in space-time $\R^{3+1}$ (the first three dimensions are interpreted as spatial dimensions and the last one is time) whose elements are given by composition of a spatial rotation $A \in \text{SO}(3)$, uniform motion with velocity $v\in \R^3$ and translations in space and time by a vector $(a,s)\in \R^{3+1}$. As a matrix group , the elements are given by
$$\begin{pmatrix} A & v & a \\ 0 & 1 & s \\ 0 & 0 & 1 \end{pmatrix}, \qquad A \in \textup{SO}(3), v,a \in \R^3, s \in \R$$
\noindent

Now let us consider the subgroup $H$ given by zero-time $s=0$. It is a closed codimension one subgroup. The pair $(G,H)$ can be useful to describe transformations that are time-preserving (or isochronous). This example serves as an inspiration for the concept of $b$-Lie group. Other groups that enter this category include the Heisenberg group. This notion will be useful also to study more sophisticated groups such as the group of gauge transformations of flat connections over a Riemann surface (see for instance \cite{atiyahbott} when the considered bundle is a $b$-bundle).

 The study of $b$-manifolds has its origins in the calculus on manifolds with boundary to give a conceptual approach to the Atiyah-Patodi-Singer index theorem in terms of the classical Atiyah-Singer theorem \cite{Melrose}. The language of $b$-tangent bundles  was also used in the extension of the deformation quantization scheme to manifolds with boundary \cite{NestandTsygan}.
$b$-Symplectic structures and normal forms for group actions on them have been intensely studied by several authors (see for instance, \cite{guimipi}, \cite{guimipi2}, \cite{gmps},\cite{marcutosorno1}, \cite{gualtierili},\cite{gualtierietal}). However, the notion of $b$-structures on Lie groups has yet to be treated. In this short article we fill this gap in the literature.

We examine $b$-structures on Lie groups and prove that they are very rigid, with few examples, due to the restrictive nature of the symmetries. A mild generalization of these groups shows up when the $b$-tangent bundle is replaced by a general Lie algebroid and the manifolds in question are $E$-manifolds \cite{evageoff} and \cite{nt}. We do not consider those here.

For $b$-Lie groups, we study the reduction of the structure of the $b$-cotangent bundle of the Lie group by the lift of the action of a subgroup $H$ identified as the critical set of the $b$-symplectic structure. These $b$-cotangent models generalize the models for integrable systems treated in \cite{km} and is a starting point to understand Poisson reduction theory in the $b$-context (see also \cite{annaevaroisin} as the $b$-analogue of the Marle-Guillemin-Sternberg normal form in \cite{Marle} and \cite{GandS}). In particular, these \emph{cotangent models} can be useful to solve some of the open problems in integrable systems both concerning normal form theory and applications to quantization (see \cite{bmmt} and \cite{mm}).

For the standard cotangent bundle of a Lie group, the Poisson geometry of $T^*B$ with $B$ a principal bundle and all the underlying theoretical framework has been used by Sternberg \cite{sternberg}, Weinstein \cite{Weinstein}, Montgomery (\cite{montgomery, montgomery2}), Marsden and Ratiu \cite{montgomeryetal}, and others to describe the mathematical formulation of particles in Yang–Mills fields. The description of the bundle picture started with Sternberg \cite{sternberg} and Weinstein \cite{Weinstein} in their investigations of the symplectic geometry of the phase space of a classical coloured particle in an external Yang-Mills field. The phase space of such a particle is a vector bundle over a standard single-particle phase space with fiber the vector space of colour charges. Though in this article we only consider the case when the principal bundle is associated to the quotient $G/H$, similar applications for different principal bundles should yield other possible applications to physics. In particular, applications for Yang-Mills fields on manifolds with boundary are envisaged.

\textbf{Acknowledgements:} We are enormously thankful to the referee for an important correction and numerous suggestions that have improved this article. In particular, their suggestions allowed us to connect with the classical theory of Yang-Mills which has been eye-opening in terms of potential future applications.

\section{Preliminaries}

A \textbf{$b$-manifold} is a pair $(M,Z)$, consisting of an oriented manifold $M$ and an oriented hypersurface $Z\subset M$.
A \textbf{$b$-vector field} on a $b$-manifold $(M,Z)$ is a vector field which is tangent to $Z$ at every point $p\in Z$.

If $f$ is a local defining function for $Z$ on some open set $U\subset M$ and $(f,z_2,\ldots,z_{n})$ is a chart on $U$, then the set of $b$-vector fields on $U$ is a free $C^\infty(U)$-module with basis
\begin{equation}\label{localframe}
\left(f \frbd{f}, \frbd{z_2},\ldots, \frbd{z_{n}}\right).
\end{equation}
\noindent
We call the vector bundle associated to this locally free $C_M^\infty$-module the \textbf{$b$-tangent bundle} and denote it $^b TM$.
We define the \textbf{$b$-cotangent bundle} $^b T^*M$ of $M$ to be the vector bundle dual to $^b TM$. The sheaf of sections of $\Lambda^k (\!\, ^b T^*M)$ is denoted $^b\Omega^k$ and its elements are called {\bf $b$-forms of degree $k$}.

The classical exterior derivative $d$ on the complex of (smooth) $k$-forms extends to the complex of $b$-forms in a natural way. Indeed, any $b$-form $\omega$ of degree $k$ can locally be written in the form
$\omega=\alpha\wedge\frac{df}{f}+\beta$ where $\alpha\in\Omega^{k-1}(M), \beta\in\Omega^k(M)$. Here $f$ is a local defining function of $Z$ and $\frac{df}{f}$ is the $b$-form of degree $1$ dual to $f \frbd{f}$ in a frame of the form \eqref{localframe}. We then define the exterior derivative
$d\omega:=d\alpha\wedge\frac{df}{f}+d\beta$ (see \cite{guimipi2} for details).
In order to have a Poincar\'{e} lemma for $b$-forms, we enlarge the set of smooth functions and consider the set of {\bf $b$-functions} $^b C^\infty(M)$, which consists of functions with values in $\R \cup \{\infty\}$ of the form $c\,\textrm{log}|f| + g$,
 where $c \in \mathbb{R}$, $f$ is a defining function for $Z$, and $g$ is a smooth function. The differential operator $d$ on this space is defined as: $d(c\,\textrm{log}|f| + g):= c \frac{df}{f} + d g$,
where $d g$ is the standard de Rham derivative.

\begin{definition}\label{bsymp}
Let $(M,Z)$ be a $b$-manifold, with $Z$ the critical hypersurface. We say a a closed $b$-form of degree $2$, $\omega\in\,\!^b\Omega^2(M)$, is \textbf{$b$-symplectic} if $\omega_p$ is of maximal rank as an element of $\Lambda^2(\,\! ^b T_p^* M)$ for all $p\in M$.
\end{definition}

A $b$-symplectic form defines a symplectic form away from $Z$, implying in particular that $M$ has even dimension. The local structure of $b$-symplectic forms is well-understood thanks to a Darboux theorem in this context:

\begin{theorem}[\textbf{$b$-Darboux theorem, \cite{guimipi2}}]\label{theorem:bDarboux}
Let $\omega$ be a $b$-symplectic form on $(M^{2n},Z)$. Let $p\in Z$. Then we can find a local coordinate chart
$$(x_1,y_1,\ldots,x_n,y_n)$$
 centered at $p$ such that the hypersurface $Z$ is locally defined by $y_1=0$ and
$$\omega=dx_1\wedge\frac{d y_1}{y_1}+\sum_{i=2}^n dx_i\wedge dy_i.$$
\end{theorem}

\subsection{\textbf{$b$-Cotangent lifts}}

Recall (see for instance, \cite{guilleminandsternberg}) that given any group action the cotangent lift of this action is Hamiltonian with respect to the canonical symplectic structure on $T^\ast M$ with moment map $\mu$:
\begin{equation*}\label{eqn:lift}
\langle\mu(p),X \rangle := \langle \lambda_p ,X^\#|_{p} \rangle =\langle p,X^\#|_{\pi(p)}\rangle,
\end{equation*}
 where $X^\#$ is the fundamental vector field of the action.
In \cite{guimipi2} it was noted that, analogous to the symplectic case, the $b$-cotangent bundle comes equipped with a canonical $b$-symplectic form.

\begin{definition}\label{canbsymp}
Let $(M,Z)$ be a $b$-manifold. Then we define a one $b$-form $\lambda$ on $^b T^*M$, considered as a $b$-manifold with critical hypersurface $^b T^*M|_Z$, in the following way:
 \begin{equation*}\label{canonical_liouville}
 \langle \lambda_p, v\rangle:= \langle p, (\pi_p)_\ast (v)\rangle, \qquad p\in \,\!^b T^*M, v \in  \,\!^b T_p(\,\!^b T^*M)
 \end{equation*}
We call $\lambda$ the {\bf $b$-Liouville form}. The negative differential
$$\omega = - d \lambda$$
is the {\bf canonical $b$-symplectic form} on $^b T^*M$.
\end{definition}

Using this form the canonical $b$-cotangent lift was defined as follows \cite{km}:

\begin{definition}\label{bcotlift}
Consider the $b$-cotangent bundle $^b T^\ast M$ endowed with the canonical $b$-symplectic structure. Assume that the action $\rho$ of $G$ on $M$ preserves the hypersurface $Z$, i.e. $\rho_g$ is a $b$-map for all $g\in G$. Then the lift of $\rho$ to an action on $^b T^\ast M$ is well-defined:
$$\hat{\rho}:G\times ^b \! T^\ast M \to \,^b \! \,T^\ast M : (g,p)\mapsto \rho^\ast_{g^{-1}}(p).$$
Moreover, it is $b$-Hamiltonian with respect to the the canonical $b$-symplectic structure on $^b T^\ast M$. We call this action together with the underlying canonical $b$-symplectic structure the {\bf canonical $b$-cotangent lift}.
\end{definition}
%
%
%

\section{Definition and motivating examples}
 In the symplectic case, reducing $T^*G$ by the action of $G$ yields the Lie-Poisson structure on $\mathfrak{g}^*$. In this paper we study the analogue in the $b$-context: In order to lift an action to $^b T^*G$ (Definition \ref{bcotlift}), we have to demand that the action leaves the critical hypersurface invariant. This motivates us to consider the setting where the critical hypersurface of the $b$-structure is a codimension one Lie subgroup $H$ and we consider the action of $H$ on $G$ by translations.

This is the content of the next definition:

\begin{definition}
A $b$-manifold $(G,H)$, where $G$ is a Lie group and $H \subset G$ is a closed codimension one subgroup \footnote{This is equivalent to $H$ being an embedded Lie subgroup.} is called a {\bf $b$-Lie group}.
\end{definition}


\begin{example}\label{e2}
The special Euclidean group of orientation-preserving isometries in the plane is the semidirect product
$$\textup{SE}(2) \cong \textup{SO}(2) \ltimes T(2)$$
\noindent
where $T(2)$ are translations in the plane. Recall that we can identify $\textup{SE}(2)$ with matrices of the form
$$\begin{pmatrix} A & b \\ 0 & 1 \end{pmatrix}, \qquad A \in \textup{SO}(2), b \in \R^2$$
\noindent
Then $T(2)$ (identified with $\{I\} \times T(2) \subset \textup{SE}(2)$) is a closed codimension 1 subgroup and the pair $(\textup{SE}(2), T(2))$ is a $b$-Lie group.
\end{example}

\begin{example}
The Galilean group $G$ is the group of transformations in space-time $\R^{3+1}$ (the first three dimensions are interpreted as spatial dimensions and the last one is time) whose elements are given by composition of a spatial rotation $A \in \text{SO}(3)$, uniform motion with velocity $v\in \R^3$ and translations in space and time by a vector $(a,s)\in \R^{3+1}$.
As a matrix group, the elements are given by
$$\begin{pmatrix} A & v & a \\ 0 & 1 & s \\ 0 & 0 & 1 \end{pmatrix}, \qquad A \in \textup{SO}(3), v,a \in \R^3, s \in \R$$
\noindent
The subgroup $H$ given by $s=0$ (which corresponds to fixing time) is a closed codimension one subgroup and hence the pair $(G,H)$ is a $b$-Lie group.
\end{example}

\begin{example}\label{e3}

We consider the $(2n+1)$-dimensional Heisenberg group $H_{2n+1}(\R)$ given by matrices of the form

$$\begin{pmatrix} 1 & a & c \\ 0 & I_n & b \\ 0 & 0 & 1 \end{pmatrix}, \qquad a\in \R^{1\times n},b \in \R^{n\times 1},c \in \R$$
\noindent
The subgroup $\Gamma$ of matrices of the form

$$\begin{pmatrix} 1 & 0 & k \\ 0 & I_n & 0 \\ 0 & 0 & 1 \end{pmatrix}, \qquad k \in \Z$$
\noindent
is central, hence normal, and so we can consider $G:=H_{2n+1}(\R)/\Gamma$. This is a well-known example of a non-matrix Lie group. Now fixing one component $a_i=0$ or $b_i =0$ yields a closed codimension one subgroup of $G$.

\end{example}

\section{Main results}
Let us consider the action of $H$ on $G$ by left translations. This action is obviously free and since $H$ is closed, it is also proper. Therefore, the right coset space $G/H$ can be given the structure of a smooth manifold such that the projection $\pi: G \to G/H$ is a smooth submersion. Moreover, it is well-known that $\pi$ turns $G$ into a principal $H$-bundle.

For future reference we summarize these facts in the following lemma:

\begin{lemma}\label{principalbundle}
Let $(G,H)$ be a $b$-Lie group. The projection $\pi: G \to G/H$ is a principal $H$-bundle; in particular $G$ is semilocally around $H$ a product
$$\pi^{-1}(V) \cong V \times H$$
\noindent
for some open neighborhood $V$ of $[e]_\sim$ in $G/H$ and under this diffeomorphism $\pi$ corresponds to the projection onto the first component.
\end{lemma}

Note that by taking a coordinate $\varphi$ on $V$ centered at $[e]_\sim$, we obtain a global defining function $\varphi \circ \pi$ for the critical hypersurface $H$.
Any local trivialization in the lemma gives rise to a natural projection $\pi_{\h}:\g \to \h$ from the Lie algebra $\g$ of $G$ to the Lie algebra $\h$ of $H$ but as we want this projection to be compatible with the bundle structure we need to consider a particular projection $\pi_{\h}:\g \to \h$ which is compatible with the bundle structure. In order to do that, we choose a projection using a connection.

Let us introduce the subbundle $\H$ of $TG$ whose fibre $\H_g$ at $g\in G$ is given by the corresponding right-shift of the Lie algebra $\mathfrak{h}$ of $H$, $\H_g = (R_g)_* \h$. Let $\pi_\H: TG \to \H$ be the natural projection onto $\H$ given by conjugating $\pi_\h$ with left-translations. Recall that $\pi: G \to G/H$ induces a surjective bundle morphism $\pi_*: TG \to T(G/H)$ and at each fibre $T_g G$ the kernel is $\H_g$.

Actually such a projection $\pi_{\mathcal{H}}$ has to be compatible with the principal bundle structure of $\pi:G\rightarrow G/H$ and we construct one using a principal $H$-connection which we describe in the subsection below.

\subsection{Basics in classical reduction theory of the cotangent bundles}

In this section we describe the general theory of the reduction theory of the cotangent bundle of a principal bundle as was developed in Montgomery's PhD thesis \cite{montgomery} and apply these classical results of this theory to the particular case when the principal $H$-bundle is given by $\pi$ : $G\rightarrow G/H.$
This paves the way to develop a reduction theory for the Poisson structure which we revise at the end of the section.
 For more details confer Section $1.1$ \cite{montgomery}.

 The reader may also wish to consult \cite{montgomery2}, \cite{sternberg} and \cite{Weinstein}.

Let us choose a principal $H$-connection $\theta\in\Omega^{1}(G;\mathfrak{h})$ for $\pi$ : $G\rightarrow G/H$, i.e. an $\mathfrak{h}$-valued 1-form on $G$ such that

$$\displaystyle \theta\left(\frac{\mathrm{d}}{\mathrm{dt}}_{\vert{t=0}} e^{tX}g\right)=X \quad \mathrm{and} \quad \theta\left(\left(L_{h}\right)_{*}v\right)=\mathrm{Ad}_{h}(\theta(v))$$
\noindent for all $g\in G, X\in \mathfrak{h}$ and $v\in TG$.

Setting
$$\pi_{\mathcal{H}}(v_{g})=(R_{g})_{*}\theta(v_{g}) \quad \textrm{for all} \quad g\in G, X\in \mathfrak{h},$$
\noindent each principal $H$-connection $\theta$ determines an Ehresmann connection \cite{connection} $\pi_{\mathcal{H}}=\pi_{\mathcal{H}}^{\theta}$ for $\pi:G\rightarrow G/H$, i.e. a projection

$$
\pi_{\mathcal{H}}:TG\rightarrow \mathcal{H}.
$$

Let us denote by $\zeta$ the infinitesimal action of $\mathfrak{h}$ on $G$ corresponding to the action of $H$ on $G$ by left multiplications:

$$\zeta:G\times \mathfrak{h}\rightarrow TG,\ (g,\ X)\mapsto\zeta_{g}^{X}:=\frac{d}{dt}\left(e^{tX}g\right)_{t=0}
$$

\noindent We can use that $\mathcal{H}$ and $\ker\theta$ are $H$-invariant subbundles of $TG$. Also, their annihilators $\mathcal{H}^{\circ}$ and $\ker{\theta}^{\circ}$ are $H$-invariant subbundles of $T^{*}G$.

The procedure described in the propositions below is often called in the literature as \emph{minimal coupling}.
 Specialising classical results from the reduction theory of the cotangent bundle to a principal bundle to the case of $\pi$ : $G\rightarrow G/H$, one gets the following proposition:

\begin{proposition}\label{prop1} Each principal $H$-connection $\theta$ for $\pi$ : $G\rightarrow G/H$, determines a diffeomorphism
$$
\phi_{\theta}\ :\ TG\rightarrow \ker{\theta}\oplus(G\times \mathfrak{h})\equiv \ker{\theta}\times \mathfrak{h}
$$
$$
v_{g}\mapsto(v_{g}-\zeta_{g}(\theta_{g}v_{g}),\ \theta_{g}v_{g})\ .
$$
Additionally, this diffeomorphism is $H$-equivariant. It intertwines:

\begin{itemize}

\item The tangent lift $H\times TG\rightarrow TG, (h,\ v_{g})\rightarrow(L_{h})_{*}v_{g}$, and

\item  the diagonal action $H\times \ker{\theta}\times \mathfrak{h}\rightarrow \ker{\theta}\times \mathfrak{h}, (h,\ (u,\ X))\mapsto((L_{h})_{*}u,\ Ad_{h}X)$ .
\end{itemize}
\end{proposition}
This result has a dual equivalent as follows:

\begin{proposition} \label{prop2}
Each principal $H$-connection $\theta$ for $\pi$ : $G\rightarrow G/H$, determines a diffeomorphism
$$
\psi_{\theta}:T^{*}G\rightarrow \mathcal{H}^{\circ}\oplus(G\times \mathfrak{h}^{*})\equiv \mathcal{H}^{\circ}\times \mathfrak{h}^{*}
$$
$$
\alpha_{g}\mapsto(\alpha_{g}-\alpha_{g}\circ\zeta_{g}\circ\theta_{g},\ \alpha_{g}\circ\zeta_{g})\ .
$$
Furthermore, this diffeomorphism is $H$-equivariant. Indeed, it intertwines:

 \begin{itemize}
   \item The cotangent lift $H\times T^{*}G\rightarrow T^{*}G, (h,\ \alpha)\rightarrow(L_{h})_{*}\alpha$, and

   \item the diagonal action $H\times \mathcal{H}^{\circ}\times \mathfrak{h}^{*}\rightarrow \mathcal{H}^{\mathrm{o}}\times \mathfrak{h}^{*}, (h,\ (\beta,\ \mu))\rightarrow((L_{h}^{-1})^{*}\beta,\ Ad_{h}^{*}\mu)$ .
 \end{itemize}

\end{proposition}

 The $H$-invariant subbundle $\mathcal{H}^{\circ}\subset T^{*}G$ has a canonical structure of principal $H$-bundle which is described in the following proposition.
\begin{proposition}\label{prop3} The annihilator of the vertical bundle $\mathcal{H}^{\circ}\subset T^{*}G$ is a principal $H$- bundle with action induced by the cotangent lifted action of $H$ on $T^{*}G$, quotient space $T^{*}(G/H)$ and bundle map $\overline{\pi}:\mathcal{H}^{0}\rightarrow T^{*}(G/H)$ defined by setting
$$
(\overline{\pi}\alpha_{g})(\pi_{*}v_{g})=\alpha_{g}v_{g},
$$
\noindent for all $g\in G, \alpha_{g}\in \mathcal{H}^{\circ}_{g}$ and $v_{g}\in T_{g}G$. Additionally, the latter fits in the following pull-back diagram which also represents a principal bundle map

 \[
\begin{tikzcd}
\mathcal{H}^{\circ}\arrow{r}{\tau_{{G}^{*}}\vert_{\mathcal{H}^\circ}} \arrow[swap]{d}{\bar{\pi}} & G \arrow{d}{\pi} \\
T^*(G/H) \arrow{r}{\tau_{G/H}^{*}} & G/H
\end{tikzcd}
\]

 Above $\tau_{G}^{*}$ : $T^{*}G\rightarrow G$ {\it and} $\tau_{G/H}^{*}$ : $T^{*}(G/H)\rightarrow G/H$ denote the standard projections.
\end{proposition}
Consequently, if the 1-form $\theta\in\Omega^{1}(G;\mathfrak{h})$ is a principal $H$-connection for $\pi$ : $G\rightarrow G/H$, then the pull-back 1-form $\tau_{G}^{*}\vert{\mathcal{H}^{\circ}}\theta\in\Omega^{1}(\mathcal{H}^{\circ};\mathfrak{h})$ is a principal $H$-connection for $\overline{\pi}:\mathcal{H}^{\circ}\rightarrow T^{*}(G/H)$ .

Now, as last preliminary result, still specializing former results of Montgomery (Section $1.1$ in \cite{montgomery}) to $\pi$ : $G\rightarrow G/H$ one gets the following description of the induced Poisson structures.

\begin{proposition}\label{prop4}

 Let $\theta$ be a principal $H$-connection for $\pi$ : $G\rightarrow G/H$. Then the induced diffeomorphism
$$
\psi_{\theta}:T^{*}G\rightarrow \mathcal{H}^{\circ}\oplus(G\times \mathfrak{h}^{*})\equiv \mathcal{H}^{\circ}\times \mathfrak{h}^{*}
$$
$$
\alpha_{g}\mapsto(\alpha_{g}-\alpha_{g}\circ\zeta_{g}\circ\theta_{g},\ \alpha_{g}\circ\zeta_{g})\ .
$$
 transforms the canonical symplectic form $\omega_{G}=-\mathrm{d}\lambda_{G}$ on $T^{*}G$ into the symplectic form on $\mathcal{H}^{\circ}\times \mathfrak{h}^{*}$ given by
$$
(\psi_{\theta})_{*}\omega_{G}=\overline{\pi}^{*}\omega_{G/H}-\mathrm{d}\lambda^{\theta}
$$
 where $\omega_{G/H}=-\mathrm{d}\lambda_{G/H}$ is the canonical symplectic form on $T^{*}(G/H)$ and $\lambda^{\theta}$ is the connection-dependent $1$-form on $\mathcal{H}^{\circ}\times \mathfrak{h}^{*}$ given by
$$
\lambda_{(\alpha,\mu)}^{\theta}(w,\ \phi)=\mu((\tau_{G}^{*}|_{\mathcal{H}^{\circ}}^{*}\theta)w)\ ,
$$
{\it for all} $\alpha\in \mathcal{H}^{\circ}, w\in T_{\alpha}\mathcal{H}^{\circ}$ {\it and} $\mu, \phi\in \mathfrak{h}^{*}$.
\end{proposition}

Let us now describe the Poisson structure on the reduced space. We do this following Section 1.1 in \cite{montgomery}.

For a general principal bundle $B$ with structural group $H$, the manifold $T^*B/H$ is a Poisson manifold. Its symplectic leaves are the spaces investigated by Weinstein in \cite{Weinstein}. The manifold $T^*B/H$ is a vector bundle over $X$ with fiber $T_x^*X\times\mathfrak{h}$ and can be identified with the coadjoint bundle $B^*\times_H \mathfrak{h}$ over $T^*X$ associated to $B^*$ (see for instance \cite{montgomery}).
We can specialize this result for the principal bundle
 $G\rightarrow G/H$.

 More concretely, the diffeomorphism in Proposition \ref{prop2} induces a diffeomorphism {$[\psi_{\theta}]$} from the reduced space $(T^*G)/H$ to the {\textbf{ coadjoint bundle}} $\mathcal{H}^{\circ}\times_H\mathfrak{h}^{*}$i.e., the total space of the vector bundle
$$
\mathcal{H}^{\circ}\times_H\mathfrak{h}^{*}:=\frac{\mathcal{H}^{\circ}\times \mathfrak{h}^{*}}{H}\rightarrow T^{*}(G/H)
$$
which is the associated bundle for the principal $H$-bundle $\overline{\pi}:\mathcal{H}^{\circ}\rightarrow T^{*}(G/H)$ and the coadjoint representation of $H$ on $\mathfrak{h}^*$.

 Next, we use the minimal coupling procedure described in Propositions \ref{prop1} and \ref{prop2} to induce a Poisson structure (as described in Proposition \ref{prop4}) on the coadjoint bundle pulling back the Poisson structure from its target.
 It was proved in \cite{montgomery} and \cite{montgomeryetal} that for general principal bundles, the symplectic leaves of this structure coincide with the ones introduced by Sternberg in \cite{sternberg}.

This can be summarized in the following diagram and theorem.

\begin{figure}[h!]
 \[
 \begin{tikzcd}
 \text{{\texttt{(Connection dependent bracket)}}} \, \mathcal{H}^{\circ}\times_H\mathfrak{h}^{*} \arrow{r}{[\psi_{\theta}]^{-1}} \arrow[swap]{dr}{p \text{ (\texttt{universal projection})}} & T^*(G)/H \, \text{\texttt{(universal bracket)}} \arrow{d}{\overline{p} \text{\texttt{(connection dependent)}}} \\
   & T^*(G/H)
 \end{tikzcd}
\]
 \caption{ { This diagram is borrowed from page 15 in \cite{montgomery} replacing a general principal bundle by $G\rightarrow G/H$. The choice of a connection determines the projection $\overline p$ through its horizontal lift (see details, for instance, in page 14 in \cite{montgomery}).} }\label{PIC}
\end{figure}
\begin{theorem}[\textbf{Poisson reduction via minimal coupling}]\label{summary} $T^*(G)/H$ has an induced Poisson structure independent of the connection inherited from the diffeomorphism of Proposition \ref{prop4} with $\mathcal{H}^{\circ}\times_H\mathfrak{h}^{*}$.
\end{theorem}

\subsection{The $H$-action on $^b T^*G$ and the $b$-minimal coupling }

Let $(G,H)$ be a $b$-Lie group and consider the action of $H$ by left translations.
In Definition \ref{bcotlift} we introduced the $b$-cotangent lift; in the present setting this is given by the following action on the $b$-cotangent bundle $^b T^* G$:
$$H \times \,\!^b T^* G \to \,\!^b T^* G:  (h, \alpha_g) \mapsto (L_{h^{-1}})^* \alpha_g.$$

Note that this action is well-defined since the left translation by $h\in H$ preserves $H$ i.e. it is a $b$-map. Moreover we define the projection $\pi_\H: \,\!^b TG \to \H$ in the analogous way by choosing a connection.

By duality we have a well defined notion of tangent lift. On the other hand,
the quotient space $(^b T^*G)/H$ is isomorphic to $\left((^b TG)/H\right)^*$ via the identification
$$(^b T^*G)/H\xrightarrow{\sim}\big((^b TG)/H\big)^* : [\alpha_g]_\sim \mapsto \big( [v_g]_\sim \mapsto \langle \alpha_g , v_g \rangle \big), \quad v_g \in \,\!^b T_g G.$$
\noindent

 So we could go straightaway to generalize the Poisson structure construction. However, let us go step by step.

Moreover, since all the maps involved are $b$-maps, one can easily reproduce all the results above replacing the cotangent bundle by the $b$-cotangent bundle.

An important caveat in this construction is that the notion of Ehresmann connection \cite{connection} has to be extended to the $b$-setting. The notion of $b$-connection is already present in the former works \cite{gmw1,gmw2} on quantization by the second author of this article.
A principal $H$-$b$-connection $\theta\in ^b \Omega^{1}(G;\mathfrak{h})$ for $\pi$ : $G\rightarrow G/H$ considering the $b$-structure of the pair $(G,H)$ is encoded on a $\mathfrak{h}$-valued 1-$b$-form on $G$ such that

$$\displaystyle \theta\left(\frac{\mathrm{d}}{\mathrm{dt}}_{\vert{t=0}} e^{tX}g\right)=X \quad \mathrm{and} \quad \theta\left(\left(L_{h}\right)_{*}v\right)=\mathrm{Ad}_{h}(\theta(v))$$
\noindent for all $g\in G, X\in \mathfrak{h}$ and $v\in ^bTG$.

This extension with the $b$-mnemonics is automatic using Melrose's $b$-calculus just by considering one forms in the complex $^b \Omega^1(G;\mathfrak{h})$.

In particular we obtain,

\begin{proposition}\label{prop1b} Each principal $H$-$b$-connection $\theta$ for the $b$-Lie group pair $(G,H)$ and the projection $\pi$ : $G\rightarrow G/H$, determines a diffeomorphism
$$
\phi_{\theta}\ :\ ^bT G\rightarrow \ker{\theta}\oplus(G\times \mathfrak{h})\equiv \ker{\theta}\times \mathfrak{h}
$$
$$
v_{g}\mapsto(v_{g}-\zeta_{g}(\theta_{g}v_{g}),\ \theta_{g}v_{g})\ .
$$
Additionally, this diffeomorphism is $H$-equivariant. Indeed, it intertwines

\begin{itemize}

\item The tangent lift $H\times ^bTG\rightarrow ^bTG, (h,\ v_{g})\rightarrow(L_{h})_{*}v_{g}$, and

\item  the diagonal action $H\times \ker{\theta}\times \mathfrak{h}\rightarrow \ker{\theta}\times \mathfrak{h}, (h,\ (u,\ X))\mapsto((L_{h})_{*}u,\ Ad_{h}X)$ .
\end{itemize}
\end{proposition}

For the $b$-cotangent lift we obtain,

\begin{proposition} \label{prop2b} Each principal $H$-$b$-connection $\theta$ for the $b$-Lie group pair $(G,H)$ and the projection $\pi$ : $G\rightarrow G/H$, determines a diffeomorphism
$$
\psi_{\theta}:^b T^{*}G\rightarrow \mathcal{H}^{\circ}\oplus(G\times \mathfrak{h}^{*})\equiv \mathcal{H}^{\circ}\times \mathfrak{h}^{*}
$$
$$
\alpha_{g}\mapsto(\alpha_{g}-\alpha_{g}\circ\zeta_{g}\circ\theta_{g},\ \alpha_{g}\circ\zeta_{g})\ .
$$
 Additionally, this diffeomorphism is $H$-equivariant intertwining:

 \begin{itemize}
   \item The cotangent lift $H\times ^b T^{*}G\rightarrow ^b T^*G, (h,\ \alpha)\rightarrow(L_{h})_{*}\alpha$, and

   \item the diagonal action $H\times \mathcal{H}^{\circ}\times \mathfrak{h}^{*}\rightarrow \mathcal{H}^{\mathrm{o}}\times \mathfrak{h}^{*}, (h,\ (\beta,\ \mu))\rightarrow((L_{h}^{-1})^{*}\beta,\ Ad_{h}^{*}\mu)$ .
 \end{itemize}

\end{proposition}

 The $H$-invariant subbundle $\mathcal{H}^{\circ}\subset ^bT^{*}G$ has a canonical structure of principal $H$-bundle which is described in the following proposition.
\begin{proposition}\label{prop3b} The annihilator of the vertical bundle $\mathcal{H}^{\circ}\subset^bT^{*}G$ is a principal $H$-bundle with action induced by the $b$-cotangent lifted action of $H$ on $^bT ^{*}G$, quotient space $^b T^{*}(G/H)$ and bundle map $\overline{\pi}:\mathcal{H}^{0}\rightarrow ^b T^{*}(G/H)$ defined by setting
$$
(\overline{\pi}\alpha_{g})(\pi_{*}v_{g})=\alpha_{g}v_{g},
$$
\noindent for all $g\in G, \alpha_{g}\in \mathcal{H}^{\circ}_{g}$ and $v_{g}\in T_{g}G$. Additionally, the latter fits in the following pull-back diagram which also represents a principal bundle map

 \[
\begin{tikzcd}
\mathcal{H}^{\circ}\arrow{r}{\tau_{{G}^{*}}\vert_{\mathcal{H}^\circ}} \arrow[swap]{d}{\bar{\pi}} & G \arrow{d}{\pi} \\
^b T^*(G/H) \arrow{r}{\tau_{G/H}^{*}} & G/H
\end{tikzcd}
\]

 Above $\tau_{G}^{*}$ : $^b T^{*}G\rightarrow G$ {\it and} $\tau_{G/H}^{*}$ : $^b T^{*}(G/H)\rightarrow G/H$ denote the standard projections.
\end{proposition}
As a consequence, if the $b$-forrm of degree 1 $\theta\in ^b \Omega^{1}(G;\mathfrak{h})$ is a principal $H$-connection for $\pi$ : $G\rightarrow G/H$, then the pull-back 1-form $\tau_{G}^{*}\vert_{\mathcal{H}^{\circ}}\theta\in ^b \Omega^{1}(\mathcal{H}^{\circ};\mathfrak{h})$ is a principal $H$-$b$-connection for $\overline{\pi}:\mathcal{H}^{\circ}\rightarrow ^b T^{*}(G/H)$ .

The following is the $b$-version of former Proposition \ref{prop4}:
\begin{proposition}\label{prop4b}

 Let $\theta$ be a principal $H$-$b$-connection for the $b$-Lie group $(G,H)$, $\pi$ : $G\rightarrow G/H$. Then the induced diffeomorphism
$$
\psi_{\theta}:^b T^{*}G\rightarrow \mathcal{H}^{\circ}\oplus(G\times \mathfrak{h}^{*})\equiv \mathcal{H}^{\circ}\times \mathfrak{h}^{*}
$$
$$
\alpha_{g}\mapsto(\alpha_{g}-\alpha_{g}\circ\zeta_{g}\circ\theta_{g},\ \alpha_{g}\circ\zeta_{g})\ .
$$
 transforms the canonical $b$-symplectic form $\omega_{G}=-\mathrm{d}\lambda_{G}$ on $^b T^{*}G$ into the $b$-symplectic form on $\mathcal{H}^{\circ}\times \mathfrak{h}^{*}$ given by
$$
(\psi_{\theta})_{*}\omega_{G}=\overline{\pi}^{*}\omega_{G/H}-\mathrm{d}\lambda^{\theta}
$$
 where $\omega_{G/H}=-\mathrm{d}\lambda_{G/H}$ is the canonical symplectic form on $^b T^{*}(G/H)$ and $\lambda^{\theta}$ is the connection-dependent $1$ $b$-form on $\mathcal{H}^{\circ}\times \mathfrak{h}^{*}$ given by
$$
\lambda_{(\alpha,\mu)}^{\theta}(w,\ \phi)=\mu((\tau_{G}^{*}|_{\mathcal{H}^{\circ}}^{*}\theta)w)\ ,
$$
{\it for all} $\alpha\in \mathcal{H}^{\circ}, w\in T_{\alpha}\mathcal{H}^{\circ}$ {\it and} $\mu, \phi\in \mathfrak{h}^{*}$.
\end{proposition}

With all these propositions set, we can now extend the classical \emph{minimal coupling} procedure to induce a Poisson structure on $(^b T^* G)/H$. We do this in the section below.


\subsection{Reduction of the canonical $b$-symplectic structure}\label{sec:red}

The cotangent bundle $T^*G$ has a canonical symplectic structure, which under the action of $G$ on itself by left translations reduces to the minus Lie-Poisson structure on $T^*G/G \cong \g^*$.
In Definition \ref{canbsymp} we have seen how to endow the $b$-cotangent bundle $^b T^*G$ with a canonical $b$-symplectic structure (with critical hypersurface $^b T^*G|_H$). What is the reduced Poisson structure on $(^b T^* G)/H$?
As a consequence of the former adaptation of the minimal coupling to the $b$-setting we obtain the $b$-analogue of the Poisson reduction theorem \ref{summary} and have a response to this natural question.

 The diffeomorphism in Proposition \ref{prop2b} induces a diffeomorphism from the reduced space $(^b T^*G)/H$ to what we could call the {\textbf{ $b$-coadjoint bundle}} $\mathcal{H}^{\circ}\times_H\mathfrak{h}^{*}$ which is the total space of the vector bundle
$$
\mathcal{H}^{\circ}\times_H\mathfrak{h}^{*}:=\frac{\mathcal{H}^{\circ}\times \mathfrak{h}^{*}}{H}\rightarrow ^b T^{*}(G/H)
$$
which is the associated bundle for the principal $H$-bundle $\overline{\pi}:\mathcal{H}^{\circ}\rightarrow ^b T^{*}(G/H)$ and the coadjoint representation of $H$ on $\mathfrak{h}^*$.
Observe that Proposition \ref{prop4b} above induces a Poisson structure via this diffeomorphism. More concretely,
 via the minimal coupling procedure described in Propositions \ref{prop1b} and \ref{prop2b} we induce a Poisson structure (as described in Proposition \ref{prop4b}) on the $b$-coadjoint bundle pulling back the Poisson structure from its target.
By applying the same argument as in the former section we obtain,

\begin{theorem}[\textbf{Poisson reduction via $b$-minimal coupling}]\label{bsummary} $^b T^*(G)/H$ has an induced Poisson structure inherited from the diffeomorphism of Proposition \ref{prop4b} with $\mathcal{H}^{\circ}\times_H\mathfrak{h}^{*}$.
\end{theorem}

It can be checked as in \cite{montgomery} and \cite{montgomeryetal} that this Poisson structure is connection-independent.
But how does this Poisson structure look like?

\textbf{Description of the Poisson structure in local coordinates}
Let $V\subset G/H$ be an open neighbourhood of $[e]_{\sim}\equiv H$ and such that $G$ trivializes as a principal $H$-bundle over $V$ (cf. Lemma \ref{principalbundle}), i.e.
$$G\supset U:=\pi^{-1}(V) \xrightarrow{\sim} H \times V$$
where the projection onto the second component corresponds to the quotient projection $\pi$; in particular the critical hypersurface $H$ gets mapped to $H \times {[e]_\sim} \subset H \times V$ and the $b$-cotangent bundle over $U$ splits in the following way:
\begin{equation*}\label{cotid}
 \,\!^b T^* U \cong T^* H \times \,\!^b T^* V.
\end{equation*}
Then the canonical $b$-symplectic structure $\omega_0$ on $^b T^* U$ is the product of the canonical symplectic structure $\omega_1$ on $T^* H$ and the canonical $b$-symplectic structure $\omega_2$ on $\,\!^b T^* V$. Denoting the Poisson tensor corresponding to $\omega_i$ by $\Pi_i$,
$$\Pi_0= \Pi_1 + \Pi_2.$$ The action of $H$ on $^b T^* U \cong T^* H \times \,\!^b T^* V$ is given by the standard cotangent lift of left translations by $H$ on $T^* H$ times the identity on $^b T^* V$. For the corresponding quotient projections $\pi_0: \,\! ^b T^*U\to( \,\!^b T^*U)/H$ and $\pi_0': T^*H \to(T^*H)/H$ we therefore have $\pi_0 = \pi_0' \times \text{id}_{^b T^* V}$. Hence the following Poisson structure is induced on $(\,\!^b T^* U)/H$ is
$$\Pi_{\text{red}} = (\pi_0)_*\Pi_0 = (\pi_0)_*(\Pi_1 + \Pi_2)=(\pi'_0)_*\Pi_1+\Pi_2.$$
Now note that $(\pi'_0)_*(\Pi_1)$ is the minus Lie Poisson structure on $\h^*$ if we identify $(T^*H)/H \cong \h^*$.

This local description coincides with the induced Poisson structure by the analogue of the Poisson reduction theorem -theorem \ref{summary}- (check the local coordinates approach in \cite{montgomery} and \cite{montgomeryetal}) proving,
\begin{theorem}[reduced $b$-Poisson structure (local expression)]
Let $^b T^*G$ be endowed with the canonical $b$-Poisson structure. Then the Poisson reduction under the cotangent lifted action of $H$ by left translations induced by the minimal coupling procedure on $(^b T^*G)/H$ is \textbf{locally equivalent} to the Poisson structure
$ (\h^* \times \,\!^b T^*(G/H),\, \Pi^-_\text{L-P} + \Pi_{b\text{-can}})$
\noindent
where $\Pi^-_\text{L-P} $ is the minus Lie-Poisson structure on $\h^*$ and $\Pi_{b\text{-can}}$ is the canonical $b$-symplectic structure on $^b T^*(G/H)$, where $G/H$ is viewed as a $b$-manifold with critical hypersurface the point $[e]_\sim$.
\end{theorem}

\begin{example}
We return to Example \ref{e2} of the special Euclidean group $\textup{SE}(2)$. Since $T(2)$ is abelian, the Lie-Poisson structure on the dual of its Lie algebra is zero. Hence $^b \,\!T^* (\textup{SE}(2))$ reduces under the action of $T(2)$ to
$$(^b \,\!T^* (\textup{SE}(2)) / T(2), \,\Pi_\text{red}) \cong (\R^2 \times \,\!^b T^* (\textup{SO}(2)), \,0 + \Pi_{b\text{-can}}),$$
\noindent

 where $\Pi_{b\text{-can}}$ is the canonical $b$-Poisson structure on $^b T^* (\textup{SO}(2))$, i.e. identifying $\textup{SO}(2)\cong \mathbb{S}^1$ in the usual way and letting $\varphi$ be the angle, $(\varphi,p)$ a $b$-canonical chart in a neighborhood of $\{\varphi=0\}$, then in these coordinates
$$\Pi_\text{red} = \varphi \frac{\partial}{\partial \varphi} \wedge \frac{\partial}{\partial p}.$$
\end{example}

\end{document}